\documentclass[11pt]{amsart}

\usepackage{amsmath,mathrsfs}

\newtheorem{Theorem}{Theorem}[section]

\newtheorem{Proposition}[Theorem]{Proposition}

\numberwithin{equation}{section}

\begin{document}

\title[Preserving Equilibria Structure of Symplectic Integrators]
 {Preservation of stability properties near fixed points of linear hamiltonian systems by symplectic integrators}

\author{Xiaohua Ding}

\address{Department of Mathematics, Harbin Institute of Technology(Wei hai),
Weihai 264209, P.R. China}

\email{mathdxh@yeah.net}

\author{ Hongyu Liu}

\address{Department of Mathematics, University of Washington, Box 354350, Seattle, WA 98195, USA}

\email{hyliu@math.washington.edu}

\author{Zaijiu Shang}

\address{Institute of Mathematics, Academy of Mathematics and System Sciences, Chinese Academy of Sciences,
 Beijing, 100080, P.R.China}

 %\email{zaijiu@amss.ac.cn}

\author{Geng Sun}

\address{Institute of Mathematics, Academy of Mathematics and System Sciences, Chinese Academy of Sciences,
 Beijing, 100080, P.R.China}

 %\email{sung@amss.ac.cn}

\author{Lingshu Wang}

\address{Institute of Mathematics, Academy of Mathematics and System Sciences, Chinese Academy of Sciences,
 Beijing, 100080, P.R.China}

%\email{lswang@amss.ac.cn}

%\thanks{}

%\subjclass{37M15, 65P10}

\keywords{symplectic Runge-Kutta method; symplectic partitioned
Runge-Kutta method; equilibria structure; composition method;
stability}\thanks{2000 \emph{Mathematical Subject Classification.} 37M15, 65P10}

\date{}

\maketitle

\begin{abstract}
Based on reasonable testing model problems, we study the
preservation by symplectic Runge-Kutta method (SRK) and symplectic
partitioned Runge-Kutta method (SPRK) of structures for fixed
points of linear Hamiltonian systems. The structure-preservation
region provides a practical criterion for choosing step-size in
symplectic computation. Examples are given to justify the
investigation.
\end{abstract}

\maketitle

\section{Introduction and preliminaries}

Consider the $n$-degree-of-freedom (d.o.f) Hamiltonian system
\begin{equation}\label{eq:H-system}
\dot{\mathbf{z}}=J\nabla H(\mathbf{z}),\qquad J=\left[
\begin{array}{cc}
  0 & I \\
  -I & 0 \\
\end{array}
\right],
\end{equation}
where $H$ is a smooth scalar function of the extended phase space
variables $\mathbf{z}\in \mathbb{R}^{2n}$, denoting the
Hamiltonian, and $J$ is the Poisson matrix with $I$ the
$n\times n$ identity matrix. By introducing the canonically
conjugate variables, $\mathbf{z}=(q, p)$, the above system can be
rewritten as
\begin{equation}\label{eq:H-system 2}
\dot{q}=\partial H/\partial p,\quad \dot{p}=-\partial H/\partial q,
\end{equation}
where $q\in \mathbb{R}^n$ represents the configuration coordinates
of the system and their canonically conjugate momenta $p\in
\mathbb{R}^n$ represents the impetus gained by movement. As is
well-known, Hamiltonian systems are introduced as a type of system
for which the existence of conservative quantities are automatic.
System~(\ref{eq:H-system 2}) possesses two remarkable properties:
\begin{enumerate}
\item[(1)] the solutions preserve the Hamiltonian, i.e.,
\begin{equation}\label{eq:preservation H}
 \frac{d H}{d t}=0;
\end{equation}

\item[(2)] the corresponding flow is symplectic, i.e.,
\begin{equation}\label{eq:preservation sym}
\frac{d}{d t} [d p\wedge d q]=0.
\end{equation}
\end{enumerate}

In the last two decades, enormous attention has been paid to
numerical methods which preserve the symplecticity, namely,
symplectic integrators for Hamiltonian systems; we refer to the monographs Hairer~\emph{et al.}~
\cite{HaiLubWan} and Sanz-Serna\&Calvo~\cite{SanCal} for details
and related literature. Theoretical
analysis together with numerous numerical experiments has shown that
symplectic integrator not only produces improved qualitative
numerical behaviors, but also allows for a more accurate long-time
scale computation than with general-purpose methods. In the
symplectic integration study, a widely recognized fact is that the
symplecticity of a numerical integrator has little to do with its
step-size. Particularly, for SRK and SPRK methods, their
symplecticities are only related to the coefficients (see Section~2
below). Therefore, in practical computations, one usually resorts to
the classical stability analyses to find a suitable range for
choosing numerical step-sizes. However, in a recent paper
\cite{HonLiuSun}, it is shown that in some cases even the step-size
of the symplectic Euler method satisfies the classical linear
stability requirements, one can still get periodic-two numerical
solutions, or even chaotic solutions. That means, we need to
require more stringent conditions on step-sizes in addition to the classical
stability considerations in simulations of Hamiltonian flows,
even for symplectic integrators. In the
present paper, we make a first step towards such investigation by
studying the influence induced by the numerical discretization on
the equilibria structure of the underlying Hamiltonian system. It is
recalled that for a general ODE of the form
\[
\dot{\mathbf{z}}=f(\mathbf{z}),\quad \mathbf{z}\in \mathbb{R}^m, \
f: \mathbb{R}^m\mapsto \mathbb{R}^m,
\]
it may admit the presence of equilibrium point, namely,
$\tilde{\mathbf{z}}\in \mathbb{R}^m$ such that
$f(\tilde{\mathbf{z}})=0$, and the eigenvalues of the corresponding
stability matrix $\nabla_\mathbf{z} f(\tilde{\mathbf{z}})$ determine
the type of the equilibrium point and its stability properties.

In the sequel, we are mainly concerned with the Runge-Kutta (RK)
methods and partitioned Runge-Kutta (PRK) methods. Henceforth, we
customarily refer to an $s$-stage RK method by the triple
$\mathcal{R}_s=(\mathnormal{A}, b, c)$, with
$\mathnormal{A}=(a_{ij})_{i,j=1}^s$, $b=(b_i)_{i=1}^s$ and
$c=(c_i)_{i=1}^s$ being, respectively, the coefficient matrix,
weights and abscissae, and an $s$-stage PRK method by the pair
$\mathcal{R}_s^{(1)}-\mathcal{R}_s^{(2)}$. Next, we would like to
review some of the classical linear stability concepts and by
tracing the origins we can thus draw forth our motivations for the
current work. The probably most well known A-stability is introduced
by Dahlquist in 1960's (see, e.g., \cite{HaiWan}). Applying
$\mathcal{R}_s$ to the famous Dahlquist test equation
\begin{equation}\label{eq:Dahlquist}
y'=\lambda y,\quad \lambda\in \mathbb{C},\ \ \Re \lambda<0,
\end{equation}
we get the following scheme
\begin{equation}\label{eq:numerical scheme}
y_{l+1}=R(z) y_l,\quad l=0,1,2,\ldots,\ \ \mbox{and}\ \ z=\lambda h,
\end{equation}
with $R(z)$ the stability function of $\mathcal{R}_s$ (see,
Chapter IV.3, \cite{HaiWan}). It is noted that the solution to
(\ref{eq:Dahlquist}) asymptotically decays to zero as $t\rightarrow
\infty$, and in order for the numerical scheme~(\ref{eq:numerical
scheme}) to yield such qualitative behavior without any restriction
on the step size $h$, we naturally require that
\begin{equation}\label{eq:A stability}
|R(z)|<1,\quad \mbox{for any $h>0$}.
\end{equation}
Methods satisfying (\ref{eq:A stability}) are called A-stable, and
this concept has been playing an indispensable role in the numerical
field. Apparently, one can derive the same conclusion ~(\ref{eq:A
stability}) for $\mathcal{R}_s$ when applying it to the following
equation
\begin{equation}\label{eq:conjugate}
y'=\bar{\lambda} y,
\end{equation}
where $\bar{\lambda}\in \mathbb{C}$ is the complex conjugate to
$\lambda$ in equation~(\ref{eq:Dahlquist}). If we set
$\lambda=\alpha+\mathrm{i} \beta$ with $\alpha, \beta\in \mathbb{R}$
and $\alpha<0$, it is easy to see that
equations~(\ref{eq:Dahlquist}) and (\ref{eq:conjugate}) are
equivalent to the following system of ODE,
\begin{equation}\label{eq:system}
\begin{cases}
&\dot{x}=\alpha x-\beta y,\\
&\dot{y}=\beta x+\alpha y.
\end{cases}
\end{equation}
System~(\ref{eq:system}) has an equilibrium point $(0,0)$ and its
corresponding stability matrix is given by
\[
\mathbf{J}=\left[%
\begin{array}{cc}
  \alpha & -\beta \\
  \beta & \alpha \\
\end{array}%
\right],
\]
which has two eigenvalues $\lambda_{1,2}=\alpha\pm \mathrm{i}\beta$.
Now, we apply $\mathcal{R}_s$ to (\ref{eq:system}) and get
\begin{equation}\label{eq:s1}
\left[%
\begin{array}{c}
  x_{l+1} \\
  y_{l+1} \\
\end{array}%
\right]=Q\left[%
\begin{array}{cc}
  R(z) & 0 \\
  0 & R(\bar{z}) \\
\end{array}%
\right]Q^{-1}\left[%
\begin{array}{c}
  x_{l} \\
  y_{l} \\
\end{array}%
\right]
\end{equation}
with
\[
z=\lambda h,\ \bar{z}=\bar\lambda h,\
\lambda=\alpha+\mathrm{i}\beta \ \mbox{and} \ Q=\frac{1}{\sqrt{2}}\left[%
\begin{array}{cc}
  \mathrm{i} & -1 \\
  1 & -\mathrm{i} \\
\end{array}%
\right].
\]
Introducing the forward difference operators
\[
\delta_t^+ x_l=\frac{x_{l+1}-x_l}{h},\qquad \delta_t^+
y_l=\frac{y_{l+1}-y_l}{h},
\]
scheme~(\ref{eq:s1}) can be reformulated into
\begin{equation}\label{eq:s2}
\delta_t^+\left[%
\begin{array}{c}
  x_{l} \\
  y_{l} \\
\end{array}%
\right]=Q\left[%
\begin{array}{cc}
  \frac{R(z)-1}{h} & 0 \\
  0 & \frac{R(\bar{z})-1}{h} \\
\end{array}%
\right]Q^{-1}\left[%
\begin{array}{c}
  x_{l} \\
  y_{l} \\
\end{array}%
\right],
\end{equation}
which is the discrete dynamical system approximating
(\ref{eq:system}). Obviously, $(0,0)$ is the fixed point
(equilibrium point) for (\ref{eq:s2}) and the corresponding
stability matrix is given by the coefficient matrix in
(\ref{eq:s2}), which is seen to possess two eigenvalues,
$\lambda_{h,1}=(R(z)-1)/h, \lambda_{h,2}=\bar{\lambda}_{h,1}$. Now,
it is readily seen that an A-stable RK method preserves the
equilibria structure unconditionally after discretization. In the
reverse, if we want an RK method to give rise to certain preservation
of the equilibria structure of the underlying ODE, we are naturally
led to the condition~(\ref{eq:A stability}). From above analyses, we
can further easily deduce that the so-called stability region for an
RK method (see, Chapter IV.3, \cite{HaiWan}) is indeed the set of
those step sizes with which the RK method can preserve the
equilibria structure of the test equation~(\ref{eq:system}).

Next, we briefly mention one more example to illustrate the close
connection between the so-called P-stability and the preservation of
equilibria structure. In the theory of orbital mechanics, many
problems are formulated as some special second order ODE
\begin{equation}\label{eq:second order}
\frac{d^2 u}{d t^2}=f(t, u),\quad u(t_0)=u_0,\quad
\dot{u}(t_0)=u_0',
\end{equation}
which is often known in advance to have periodic solution. In
\cite{LamWat}, Lambert and Watson introduced the concept of
P-stability for numerical methods to solve equation (\ref{eq:second
order}). It is remarked that in \cite{LamWat}, the P-stability is
defined for multistep methods, while it is equally applicable to RK
method. The test model problem given in \cite{LamWat} is
\begin{equation}\label{eq:test 2nd order}
\ddot{u}=-\lambda^2 u,\quad \lambda\in \mathbb{R}.
\end{equation}
Set $\dot{u}=v$, then equation~(\ref{eq:test 2nd order}) can be
reformulated into a first order system
\begin{equation}\label{eq:2nd system}
\begin{cases}
& \dot{u}=v,\\
& \dot{v}=-\lambda^2 u.
\end{cases}
\end{equation}
In like manner as above for A-stability, one can show that the
P-stable RK methods for (\ref{eq:test 2nd order}) preserve the
equilibria structure of (\ref{eq:2nd system}) unconditionally and
vice versa. However, we will not explore more details on this aspect
and now turn our study to Hamiltonian system ~(\ref{eq:H-system}).

The equilibrium points of system~(\ref{eq:H-system}) are those
$(\tilde p, \tilde q)\in \mathbb{R}^n\times \mathbb{R}^n$ such that
$\nabla H(\tilde p, \tilde q)=0$ and the corresponding stability
matrices are given by
\[
\mathbf{J}_S (\tilde p, \tilde q):=\frac{\partial(-\partial
H/\partial q,
\partial H/\partial p)}{\partial (p, q)}(\tilde p, \tilde q).
\]
To ease our study, we start with linear Hamiltonian systems, and
then propose several ways to extend our analyses to
nonlinear case. We first designate some model
problems for our investigation. Apparently, these model problems
should feature the Hamiltonian systems. In \cite{Mos}, it is shown
that an $n$ d.o.f linear Hamiltonian system can be canonically
transformed into a Hamiltonian system consisting of $n$ 1-d.o.f
subsystems and these subsystems assume the following standard forms,
\begin{equation}\label{eq:a}
\begin{cases}
\dot{p}=& -\beta^2 q,\\
\dot{q}=& p,\qquad\quad \beta>0
\end{cases}
\end{equation}
or
\begin{equation}\label{eq:b}
\begin{cases}
\dot{p}=& -\beta p,\\
\dot{q}=& \beta q,\qquad\quad \beta>0.
\end{cases}
\end{equation}
Therefore, it is natural for us to take systems~(\ref{eq:a}) and
(\ref{eq:b}) as our testing model problems and which will be used
throughout. Clearly, $(0, 0)$ is their equilibrium point, and for
system~(\ref{eq:a})
\[
\mathbf{J}_s=\left[
\begin{array}{cc}
  0 & -\beta^2 \\
  1 & 0 \\
\end{array}
\right]
\]
with corresponding eigenvalues $\lambda_{1,2}=\pm \mathrm{i} \beta$,
which are of elliptic type, and we call system~(\ref{eq:a}) have
elliptic equilibria structure; while for system~(\ref{eq:b})
\[
\mathbf{J}_s=\left[
\begin{array}{cc}
  -\beta & 0 \\
  0 & \beta \\
\end{array}
\right]
\]
with the eigenvalues $\lambda_{1,2}=\pm \beta$, which are of
hyperbolic type, and we call system~(\ref{eq:b}) having hyperbolic
equilibria structure. Next, we shall investigate under what
conditions, the equilibria structure can be inherited by the
corresponding numerical schemes. And what would be caused by the
numerical discretization if the structure-preservation conditions
are ruined.

In the section followed, we present some necessary and sufficient
conditions for an SRK/SPRK method and its corresponding symmetric
composition method to preserve the equilibria structures of the
testing problems. The fundamental idea is that we regard the
numerical scheme as a discrete dynamic which is an approximation to
the Hamiltonian system and therefore should have the same equilibria
structure as its continuous counterpart.  Results show that in some
cases, SRK method can preserve the equilibria structure
unconditionally, but without exception, PRK method always gives a
structure-preservation region on the positive real line for the
numerical step-size. Furthermore, it is found that if the order of
accuracy of a SPRK method is increased by symmetric composition, the
structure-preservation region will be decreased accordingly.

In section~3, we extend the analyses to nonlinear Hamiltonian
systems. Besides, we give an example which justifies the necessities
of such equilibria structure preservation study.

\section{Equilibria structure preservation of SRK/SPRK method}

In this paper, we confine our study to symplectic Runge-Kutta (SRK)
methods and symplectic partitioned Runge-Kutta (SPRK) methods. The
symplecticity conditions for $\mathcal{R}_s$ and
$\mathcal{R}_s^{(1)}-\mathcal{R}_s^{(2)}$ are, respectively, given
by (see \cite{HaiLubWan},\cite{SanCal})
\begin{align}
\mathnormal{B} \mathnormal{A}+\mathnormal{A}^T \mathnormal{B}- b
b^T=& 0, \qquad \mathnormal{B}=\mbox{diag} [b]; \\
\mathnormal{B}^{(2)} \mathnormal{A}^{(1)}+{\mathnormal{A}^{(2)}}^T
\mathnormal{B}^{(1)}- b^{(2)} {b^{(1)}}^T=& 0, \quad
\mathnormal{B}^{(i)}=\mbox{diag} [b^{(i)}]\,\,\,
(i=1,2),\\
&\hspace*{-3cm}b^{(1)}=b^{(2)}.
\end{align}

\subsection{Preserving the elliptic equilibria structure}

We first apply $\mathcal{R}_s$ to system~(\ref{eq:a}) and get the
following scheme
\begin{equation}\label{eq:scheme a RK}
\begin{cases}
&P=p_l e_s-h \beta^2 A Q, \qquad p_{l+1}=p_l-h \beta^2 b^T Q, \\
&Q= q_l e_s + h A P,\qquad \quad \, q_{l+1}=q_l+h b^T P,
\end{cases}
\end{equation}
where $P=[P_1, P_2,\ldots, P_s]^T$ and $Q=[Q_1, Q_2, \ldots, Q_s]^T$
are internal stage values, $e_s=[1,\ldots,1]^T$ and $p_l\approx p(l
h), q_l\approx q(lh)$. It can be computed straightforwardly that
\[
\begin{cases}
& P=(I+h^2\beta^2 A^2)^{-1}(e_s p_l -h\beta^2 A e_s q_l ),\\
& Q=(I+h^2 \beta^2 A^2)^{-1}(e_s q_l+ h A e_s p_l),
\end{cases}
\]
and thus scheme~(\ref{eq:scheme a RK}) is equivalently reformulated
as
\begin{equation}\label{eq:discrete dynam RK for a}
\begin{cases}
& \delta_t^+ p_l=-h \beta^2 D A e_s p_l- \beta^2 D e_s q_l,\\
& \delta_t^+ q_l= D e_s p_l- h\beta^2 D A e_s q_l,
\end{cases}
\end{equation}
where $D=b^T(I+ h^2\beta^2 A^2)^{-1}$ and the forward difference
operators are defined as
\[
\delta_t^+ p_l=\frac{p_{l+1}-p_l}{h},\qquad \delta_t^+
q_l=\frac{q_{l+1}-q_l}{h}.
\]
(\ref{eq:discrete dynam RK for a}) is a discrete dynamical system
which approximates (\ref{eq:a}) and obviously, $(0,0)$ is its fixed
point (equilibrium point). The stability matrix of (\ref{eq:discrete
dynam RK for a}) is
\[
\mathbf{J}_N= \left[
\begin{array}{cc}
  -h \beta^2 D A e_s & -\beta^2 D e_s \\
  D e_s & -h \beta^2 D A e_s \\
\end{array}
\right]
\]
and its eigenvalues are
\[
\widetilde{\lambda}_{1,2}=- h \beta^2 D A e_s \pm \mathrm{i}(\beta D
e_s),
\]
which are of elliptic type. Hence, the SRK method preserves the
elliptic equilibria structure of system~(\ref{eq:a})
unconditionally, but here we note that there is a little shift of
$\widetilde{\lambda}_{1,2}$ from $\lambda_{1,2}$ due to the
numerical discretization. Next, we apply an $s$-stage SPRK method
$\mathcal{R}_s^{(1)}-\mathcal{R}_s^{(2)}$ to system~(\ref{eq:a}) and
get the following scheme
\begin{equation}\label{eq:scheme a PRK}
\begin{cases}
&P=p_l e_s-h \beta^2 A^{(1)} Q, \qquad p_{l+1}=p_l-h \beta^2 {b^{(1)}}^T Q, \\
&Q= q_l e_s + h A^{(2)} P,\qquad \quad \, q_{l+1}=q_l+h {b^{(2)}}^T
P.
\end{cases}
\end{equation}
The discrete dynamical system equivalent to (\ref{eq:scheme a
PRK}) is given by
\begin{equation}\label{eq:discrete dynam PRK for a}
\begin{cases}
& \delta_t^+ p_l=-h \beta^2 D^{(1)} A^{(2)} e_s p_l- \beta^2 D^{(1)} e_s q_l,\\
& \delta_t^+ q_l= D^{(2)} e_s p_l- h\beta^2 D^{(2)} A^{(1)} e_s q_l,
\end{cases}
\end{equation}
where
\begin{equation}\label{eq:D 1 and 2}
D^{(1)}={b^{(1)}}^T(I+ h^2\beta^2 A^{(2)}A^{(1)})^{-1}, \quad
D^{(2)}={b^{(2)}}^T(I+ h^2\beta^2 A^{(1)}A^{(2)})^{-1}.
\end{equation}
The stability matrix for (\ref{eq:discrete dynam PRK for a}) is
\[
\mathbf{J}_N= \left[
\begin{array}{cc}
  -h \beta^2 D^{(1)} A^{(2)} e_s & -\beta^2 D^{(1)} e_s \\
  D^{(2)} e_s & -h \beta^2 D^{(2)} A^{(1)} e_s \\
\end{array}
\right],
\]
and its eigenvalues are
\begin{align*}
\widetilde{\lambda}_{1,2}=&\frac{1}{2}\bigg \{-h\beta^2[D^{(1)}
A^{(2)}e_s+D^{(1)} A^{(2)}e_s]\\
& \pm \sqrt{h^2 \beta^4 [D^{(1)} A^{(2)}e_s-D^{(2)} A^{(1)} e_s]^2-4
\beta^2 [(D^{(1)}e_s)(D^{(2)}e_s)]} \bigg \}.
\end{align*}
Therefore, if the equilibria point $(0, 0)$ for (\ref{eq:discrete dynam
PRK for a}) is of elliptic type we need to require
\begin{equation}\label{eq:cond for PRK elliptic}
h^2 \beta^2 [D^{(1)} A^{(2)} e_s-D^{(2)} A^{(1)} e_s]^2-4
[(D^{(1)}e_s)(D^{(2)}e_s)]<0.
\end{equation}
Next, we introduce
\begin{equation}\label{eq:a1 and a2}
\begin{split}
a_1:=&\frac{\det (I+z^2 A^{(2)}(A^{(1)}-e_s {b^{(1)}}^T))}{\det
(I+z^2 A^{(2)}A^{(1)})},\\
a_2:=&\frac{\det (I+z^2 A^{(1)}(A^{(2)}-e_s {b^{(2)}}^T))}{\det
(I+z^2 A^{(1)}A^{(2)})},
\end{split}
\end{equation}
where $z=\beta h$, which shall also be used throughout, then
condition~(\ref{eq:cond for PRK elliptic}) can be reformulated as
\begin{equation}\label{eq:relat}
|a_1+a_2|<2.
\end{equation}
In fact, from (\ref{eq:discrete dynam PRK for a}), one can get
\begin{equation}\label{eq:relation}
\begin{split}
p_{l+1}=& (1-z^2 D^{(1)} A^{(2)}e_s) p_l-\beta z D^{(1)} e_s q_l,\\
q_{l+1}=& h D^{(2)} e_s p_l+(1-z^2 D^{(2)} A^{(1)} e_s) q_l.
\end{split}
\end{equation}
Since the method is symplectic, it has
\[
d p_{l+1}\wedge d q_{l+1}=d p_l \wedge d q_l.
\]
Substituting the relations (\ref{eq:relation}) into the above
equality, then through straightforward calculations, we obtain
\[
(D^{(1)} e_s)(D^{(2)} e_s)=[D^{(1)} A^{(2)} e_s+D^{(2)} A^{(1)}
e_s]- z^2 (D^{(1)}A^{(2)} e_s)(D^{(2)}A^{(1)} e_s),
\]
which is then substituted into (\ref{eq:cond for PRK elliptic})
to yield
\begin{equation}\label{eq:relation3}
|1-z^2 D^{(1)} A^{(2)} e_s+1-z^2 D^{(2)} A^{(1)} e_s|<2.
\end{equation}
Next, by observing that
\begin{align*}
1-z^2 D^{(1)} A^{(2)} e_s=& 1-z^2 {b^{(1)}}^T(I+ h^2\beta^2
A^{(2)}A^{(1)})^{-1} A^{(2)} e_s=a_1,\\
1-z^2 D^{(2)} A^{(1)} e_s=& 1-z^2 {b^{(2)}}^T(I+ h^2\beta^2
A^{(1)}A^{(2)})^{-1} A^{(1)} e_s=a_2,
\end{align*}
we finally arrive at the equivalence of (\ref{eq:cond for PRK elliptic})
and (\ref{eq:relat})

In summary, we have
\begin{Proposition}\label{prop:1}
For system~(\ref{eq:a}), an SRK method $\mathcal{R}_s$ preserves its
elliptic equilibria structure unconditionally, whereas for an SPRK
method, $\mathcal{R}_s^{(1)}-\mathcal{R}_s^{(2)}$, its
structure-preservation region is given by
\begin{equation}\label{eq:region1}
\{h>0: |a_1+a_2|<2\},
\end{equation}
where $a_1$ and $a_2$ are defined by (\ref{eq:a1 and a2}).
\end{Proposition}

Consider the symplectic Euler method, where $A^{(1)}=0, A^{(2)}=1$
and $b^{(1)}=b^{(2)}=1$. By (\ref{eq:a1 and a2}), it is computed
that $a_1=1-z^2, a_2=1$, and hence for symplectic Euler method to
preserve the elliptic structure of system~(\ref{eq:a}) we must
have from~(\ref{eq:region1}) that $|2-z^2|<2$, i.e., $0<\beta
h<2$.

Using the W-transformation due to Hairer and Wanner ( see Chapter~
IV.5, \cite{HaiWan}), we have
\begin{equation}\label{eq:w-transformation}
W^{-1} A^{(1)} W=X_{A^{(1)}},\qquad W^{-1} A^{(2)} W=X_{A^{(2)}},
\end{equation}
where $W$ is the generalized Vandermonde matrix, and $X_{A^{(i)}}
(i=1,2)$ are matrices possessing some standard form. In particular,
for a class of important symplectic method, Lobatto IIIA-Lobatto
IIIB pair (see \cite{sun1,sun2}),
\begin{align}
X_{A^{(1)}}=& \left(%
\begin{array}{ccccc}
  \frac 1 2 & -\xi_1 &  &  &  \\
  \xi_1 & 0 & \ddots &  &  \\
   & \ddots & \ddots & -\xi_{s-2} &  \\
   &  & \xi_{s-2} & 0 & 0 \\
   &  &  & \xi_{s-1} & 0 \\
\end{array}%
\right),\label{eq:Lobatto A}\\
X_{A^{(2)}}=& \left(%
\begin{array}{ccccc}
  \frac 1 2 & -\xi_1 &  &  &  \\
  \xi_1 & 0 & \ddots &  &  \\
   & \ddots & \ddots & -\xi_{s-2} &  \\
   &  & \xi_{s-2} & 0 & -\xi_{s-1} \\
   &  &  & 0 & 0 \\
\end{array}%
\right),\label{eq:Lobatto B}
\end{align}
where $\xi_k=1/2\sqrt{4 k^2-1}, k=1,2,\ldots,s$. In term of
W-transformation, i.e., using (\ref{eq:w-transformation}), $a_i
(i=1,2)$ in (\ref{eq:a1 and a2}) can be reformulated as
\begin{align}
a_1=& \frac{\det (I+z^2 X_{A^{(2)}}(X_{A^{(1)}}-e^{1}
{e^{1}}^T))}{\det (I+z^2 X_{A^{2}}X_{A^{(1)}})},\\
a_2=& \frac{\det (I+z^2 X_{A^{(1)}}(X_{A^{(2)}}-e^{1}
{e^{1}}^T))}{\det (I+z^2 X_{A^{1}}X_{A^{(2)}})},
\end{align}
where $e^1=[1,0,\ldots,0,1]^T$. For Lobatto IIIA-Lobatto IIIB pair,
with $X_{A^{(i)}} (i=1,2)$ given in (\ref{eq:Lobatto A}) and
(\ref{eq:Lobatto B}), it is easily verified that $a_1=a_2$ and
therefore, for Lobatto IIIA-Lobatto IIIB method to preserve the
elliptic structure of system~(\ref{eq:a}) we must have from
(\ref{eq:region1}) that $|a_1|<1$. In the following, we list some of
the computed results for Lobatto IIIA-Lobatto IIIB method in term of
its stage $s$; see Table~\ref{tab:1}.
\begin{table}[h]
  \centering
  \caption{Structure-preservation region for Lobatto IIIA-Lobatto IIIB pair}\label{tab:1}
\begin{tabular}{ccc}
  \hline
  $s$ & $a_1$ & $z$ \\
  \hline
  2 & $1-\frac 1 2 z^2$ & $0<z<2$ \\
  3 & $\frac{1-\frac{11}{24}z^2+\frac{1}{48}z^4}{1+\frac{1}{24}z^2}$ & $0<z<\sqrt{8}(\approx 2.828)$ \\
  4 & $\frac{1-\frac{7}{15}z^2+\frac{23}{900}z^4-\frac{1}{3600}z^6}{1+\frac{1}{30}z^2+\frac{1}{1800}z^4}$ & $0<z<\sqrt{42-6\sqrt{29}}(\approx 3.1127)$ \\
  \hline
\end{tabular}
\end{table}
Moreover, for $s=5$, the structure-preservation region is
$0<z<3.140328$, and for $s=10$, it is $0<z<3.141590$. We are
naturally led to the conjecture that as $s\rightarrow \infty$, the
elliptic structure-preservation region for an $s$-stage Lobatto
IIIA-Lobatto IIIB method is given by $0<z<\pi$.

Next, we consider the composition of a given basic one-step method
with different step sizes,
\begin{equation}\label{eq:symm composition}
\Psi_h=\Phi_{\gamma_s h}\circ\cdots\circ \Phi_{\gamma_1 h},
\end{equation}
where $\Phi_h$ is the basic SRK method and $\gamma_i h
(i=1,\ldots,s)$ are the composition step length (cf.
\cite{McL},\cite{Yos}). The composition is required to be symmetric,
i.e., $\gamma_i=\gamma_{s+1-i} (i=1,2,\ldots,[s/2])$, and hence
$\Psi_h$ is still an SRK method. By Proposition~\ref{prop:1}, we know
that the symmetric composition of an SRK method still preserves the
elliptic structure of system~(\ref{eq:a}) unconditionally. Next,
consider the structure-preservation region of a symmetric
composition of a SPRK method, which is known to be again a SPRK method.
To ease our study, we only take the 2-stage Lobatto IIIA-Lobatto
IIIB pair as an example, whose 4th-order symmetric composition is
given in \cite{LiuSun}, and for which we compute
\[
a_1=1-\frac 1 2 z^2+\frac 1 4
\gamma_1^2\gamma_2^2(5\gamma_1+4\gamma_2)z^4-\frac 1 8
\gamma_1^3\gamma_2(\gamma_1+\gamma_2)^2 z^6,
\]
where $\gamma_1=1/(2-2^{1/3})$ and $\gamma_2=2^{1/3}/(2-2^{1/3})$.
By $|a_1|<1$, we get
\[
0<z<\sqrt{2.48}(\approx 1.5748).
\]
Similarly, we further computed that the structure-preservation
region for the corresponding 6th-order composition method is given
by $0<z<1.1034$. Noting by Table~\ref{tab:1}, the
structure-preservation region for the underlying basic method is
$0<z<2$. Base on this example, we conjecture that as the order of a SPRK
method is increased by symmetric composition, its
structure-preservation region will be decreased accordingly. A
stringent proof of such conjecture is fraught with difficulties and is
beyond the scope of this paper.

\subsection{Preserving the hyperbolic equilibria structure}
The scheme of an SRK method $\mathcal{R}_s$ for system~(\ref{eq:b})
is read as
\begin{equation}\label{eq:discrete dynam RK for b}
\begin{cases}
& \delta_t^+ p_l=\frac{R(-z)-1}{h} p_l,\\
& \delta_t^+ q_l= \frac{R(z)-1}{h} q_l,
\end{cases}
\end{equation}
where $R(z)=\det(I+zA)/\det (I-zA)$ is the stability matrix for
$\mathcal{R}_s$. It is readily observed that for the discrete
dynamical system~(\ref{eq:discrete dynam RK for b}) to preserve the
equilibria structure of system~(\ref{eq:b}), we must have
\begin{equation}\label{eq:region for RK b}
R(-z)-1<0,\qquad \mbox{and} \qquad R(z)-1>0,
\end{equation}
i.e., the hyperbolic structure-preservation region for an SRK method
$\mathcal{R}_s$ is
\begin{equation}\label{eq:region for RK b 2}
\{h>0: \frac{\det (I-zA)}{\det (I+z A)}<1 \ \ \mbox{and} \ \
\frac{\det (I+zA)}{\det (I-zA)}>1, z=\beta h\}.
\end{equation}
Since we always have $R(-z)<1 (z=\beta h>0)$ for an A-stable RK
method. Hence, if the RK method $\mathcal{R}_s$ is A-stable, in order
to preserve the hyperbolic structure, we only need to require
\begin{equation}\label{eq:t1}
\frac{\det (I+z A)}{\det (I- z A)}>1.
\end{equation}
Now, we show that the above inequality can be equivalently reduced
to
\begin{equation}\label{eq:t2}
\det (I-z A)>0.
\end{equation}
Indeed, since $\mathcal{R}_s$ is A-stable, the eigenvalues of its
coefficient matrix must lie on the right half plane, which implies
$\det (I+z A)>0$ for $z>0$, and this together with (\ref{eq:t1})
further implies that $\det (I-zA)>0$. Next, using again the fact
that the eigenvalues of $A$ lie on the right half plane, we have
that
\[
\frac{\det (I+zA)}{\det (I-z A)}>1,
\]
if $z>0$ satisfying $\det (I-z A)>0$. Therefore, the hyperbolic
equilibria structure-preservation region of an A-stable SRK method
$\mathcal{R}_s$ is given by
\begin{equation}\label{eq:region for A-stable RK}
\{h>0: \det (I-z A)>0, \ z=\beta h\}.
\end{equation}
For example, the hyperbolic structure-preservation region for the
well-known midpoint formula, where $A=1/2$ is given by
\begin{equation}\label{eq:temp 1}
\{h>0: 1-\frac 1 2 \beta h>0\}=\{h>0: \beta h <2  \}.
\end{equation}
It is readily seen that even for SRK method which possesses good
classical stability properties, we should still require some
restrictions on its step-sizes for practical computations.
Therefore, the investigations on the equilibria structure
preservation provide a novel and useful criteria to choose step-sizes for
symplectic integrators in addition to
the classical linear stability requirements.

We now apply an $s$-stage SPRK method
$\mathcal{R}_s^{(1)}-\mathcal{R}_s^{(2)}$ to system~(\ref{eq:b}) and
get the following scheme
\begin{equation}\label{eq:discrete dynam PRK for b}
\begin{cases}
& \delta_t^+ p_l=\frac{R^{(1)}(-z)-1}{h} p_l,\\
& \delta_t^+ q_l= \frac{R^{(2)}(z)-1}{h} q_l,
\end{cases}
\end{equation}
where $R^{(i)}(z)=\det(I-zA^{(i)}+z e_s {b^{(i)}}^T)/\det
(I-zA^{(i)}) (i=1,2),$ are the stability matrices for
$\mathcal{R}_s^{(i)} (i=1,2)$. Consequently, the hyperbolic
structure-preservation region for
$\mathcal{R}_s^{(1)}-\mathcal{R}_s^{(2)}$ is
\begin{equation}\label{eq:region for RK b 2}
\{h>0: R^{(1)}(-z)<1 \ \ \mbox{and} \ \ R^{(2)}(z)>1, z=\beta h\}.
\end{equation}
For the well-known Lobatto IIIA-Lobatto IIIB pair,
$\mathcal{R}_s^{(1)}-\mathcal{R}_s^{(2)}$, since both Lobatto IIIA
method ($\mathcal{R}_s^{(1)}$) and Lobatto IIIB method
($\mathcal{R}_s^{(2)}$) are symmetric A-stable RK method, and
therefore $R^{(i)}(z)=\det(I+zA^{(i)})/\det (I-zA^{(i)}) (i=1,2).$
Similar to (\ref{eq:region for A-stable RK}),
condition~(\ref{eq:region for RK b 2}) is reduced to
\begin{equation}\label{eq:region for A-stable PRK}
\{h>0: \det (I-z A^{(i)})>0, \ i=1,2, \ z=\beta h\}.
\end{equation}
for Lobatto IIIA-Lobatto IIIB method, where in particular, we note that
$\det (I-z A^{(i)})=\det (I-z X_{A^{(1)}})=\det (I- z X_{A^{(2)}})
(i=1,2)$ with $X_{A^{(i)}}(i=1,2)$ given in (\ref{eq:Lobatto A}) and
(\ref{eq:Lobatto B}).

In summary, we have
\begin{Proposition}\label{prop:2}
The hyperbolic equilibria structure-preservation region of an SRK
method $\mathcal{R}_s$ for system~(\ref{eq:b}) is
\[
\{h>0: \frac{\det (I-zA)}{\det (I+z A)}<1 \ \ \mbox{and} \ \
\frac{\det (I+zA)}{\det (I-zA)}>1, z=\beta h\},
\]
and this condition is reduced to
\[
\{h>0: \det (I-z A)>0, \ z=\beta h\},
\]
if $\mathcal{R}_s$ is A-stable. For an SPRK method,
$\mathcal{R}_s^{(1)}-\mathcal{R}_s^{(2)}$, the hyperbolic
structure-preservation region is given by
\[
\{h>0: R^{(1)}(-z)<1 \ \ \mbox{and} \ \ R^{(2)}(z)>1, z=\beta h\},
\]
where $R^{(1)}$ and $R^{(2)}$ are, respectively, the stability
functions for $\mathcal{R}_s^{(1)}$ and $\mathcal{R}_s^{(2)}$.
\end{Proposition}

Now, we consider the structure-preservation of the symmetric
composition method~(\ref{eq:symm composition}). The 4th order
symmetric composition of mid-point formula for system~(\ref{eq:b})
is
\begin{equation}\label{eq:hyper comp}
\begin{cases}
p_{l+1}=&\displaystyle{\frac{(1-\gamma_1 z)(1-\gamma_2 z)(1-\gamma_3
z)}{(1+\gamma_1 z)(1+\gamma_2 z)(1+\gamma_3 z)}} p_l,\\
 q_{l+1}=& \displaystyle{\frac{(1+\gamma_1
z)(1+\gamma_2 z)(1+\gamma_3 z)}{(1-\gamma_1 z)(1-\gamma_2
z)(1-\gamma_3 z)}} q_l
\end{cases}
\end{equation}
with $z=\beta h$, $\gamma_1=1/(2-2^{1/3})$ and
$\gamma_2=2^{1/3}/(2-2^{1/3})$. Therefore, for scheme to preserve
the hyperbolic structure, we have
\[
1-\frac{z}{2-2^{1/3}}>0\qquad \mbox{and}\qquad
1-\frac{2^{1/3}z}{2-2^{1/3}}>0,
\]
which yields
\begin{equation}\label{eq:temp 2}
0<z<\frac{2-2^{1/3}}{2^{1/3}}(\approx 0.587401).
\end{equation}
Clearly, in comparison with (\ref{eq:temp 1}), the
structure-preservation region becomes smaller after composing. This
result applies equally to the SPRK method. Consider the following
symplectic Euler method for system~(\ref{eq:b}),
\begin{equation}\label{eq:temp}
\begin{cases}
p_{l+1}=& (1-z) p_l,\\
q_{l+1}=& \displaystyle{\frac{1}{1-z}} q_l,\qquad z=\beta h,
\end{cases}
\end{equation}
whose structure-preservation region can be verified to be
\begin{equation}\label{eq:temp 3}
0<z<1.
\end{equation}
Since the 2nd order symmetric composition of symplectic Euler method
is the midpoint formula, hence the structure-preservation region of
the 4th order symmetric composition of scheme~(\ref{eq:temp}) is
given by (\ref{eq:temp 2}), which is obviously decreased in
comparison with (\ref{eq:temp 3}).

\section{Extensions to nonlinear Hamiltonian systems and some applications}

In this section, we extend the analyses in the previous sections to
the nonlinear 1-d.o.f.
Hamiltonian system
\begin{equation}\label{eq:sect 3 1}
\begin{cases}
&\dot{p}=-{\partial H}/{\partial q}:=f(p, q),\\
&\dot{q}= {\partial H}/{\partial p}:=g(p, q),
\end{cases}
\end{equation}
with $(p, q)\in \mathbb{R}^2$. We denote by $\mathscr{E}_H$ the set
of equilibrium points for system~(\ref{eq:sect 3 1}), i.e.,
\[
\mathscr{E}_H:=\bigg\{(\tilde{p}, \tilde{q})\in \mathbb{R}^2;\
\nabla_{(p, q)} H(\tilde{p}, \tilde{q})=0\bigg\}.
\]
The stability matrix for system~(\ref{eq:sect 3 1}) is locally defined for every $(\tilde p,\tilde q)\in
\mathscr{E}_H$ as
\begin{equation}
\mathbf{J}_S(\tilde p,\tilde q)=\left[
\begin{array}{cc}
  {\partial f}/{\partial p} & {\partial f}/{\partial q} \\
  {\partial g}/{\partial p} & {\partial g}/{\partial q} \\
\end{array}%
\right](\tilde p,\tilde q),
\end{equation}
whose eigenvalues are given by
\[
\lambda_{1,2}=\pm \sqrt{D}
\]
with
\[
D:=D(\tilde p,\tilde q)=[\frac{\partial f}{\partial q}\frac{\partial g}{\partial
p}-\frac{\partial f}{\partial p}\frac{\partial g}{\partial q}](\tilde p,\tilde q).
\]
Therefore, we can make the following classification:
\begin{enumerate}
\item[(i)] If $D(\tilde{p},\tilde{q})<0$ for all $(\tilde{p},\tilde{q})\in
\mathscr{E}_H$, then the equilibrium points for
system~(\ref{eq:sect 3 1}) are all of elliptic type, and we then
call system~(\ref{eq:sect 3 1}) having elliptic equilibria
structure;

\item[(ii)] If $D(\tilde{p},\tilde{q})>0$ for all $(\tilde{p},\tilde{q})\in
\mathscr{E}_H$, then the equilibrium points for
system~(\ref{eq:sect 3 1}) are all of hyperbolic type, and we call
system~(\ref{eq:sect 3 1}) having hyperbolic equilibria structure;

\item[(iii)] If $D(\tilde{p},\tilde{q})\geq 0$ for some $(\tilde{p},\tilde{q})\in
\mathscr{E}_H$, while $D(\tilde{p},\tilde{q})\leq 0$ for some
$(\tilde{p},\tilde{q})\in \mathscr{E}_H$, we call
system~(\ref{eq:sect 3 1}) having mixed-type equilibria structure.
\end{enumerate}

As an example, we suppose that system~(\ref{eq:sect 3 1}) have hyperbolic structure
and apply an $s$-stage SRK method $\mathcal{R}_s$ for its
discretization to yield
\begin{equation}\label{eq:sect 3 2}
\begin{cases}
& P=p_l e_s+ h A F(P,Q),\qquad p_{l+1}=p_l+h b^T F(P, Q),\\
& Q=q_l e_s+ h A G(P,Q),\qquad q_{l+1}=q_l+h b^T G(P,Q),
\end{cases}
\end{equation}
where $F(P,Q)=[f(P_1,Q_1),\ldots,f(P_s,Q_s)]^T$ and
$G(P,Q)=[g(P_1,Q_1),\ldots,$ $g(P_s,Q_s)]^T$ are the internal stage
vectors. The corresponding discrete dynamical system for scheme
(\ref{eq:sect 3 2}) is
\begin{equation}\label{eq:sect 3 3}
\begin{cases}
& \delta^+ p_l= b^T F(P,Q),\\
& \delta^+ q_l= b^T G(P,Q).
\end{cases}
\end{equation}
Clearly, the points in $\mathscr{E}_H$ are still the equilibrium
points (or fixed points) of (\ref{eq:sect 3 3}). The stability
matrix for (\ref{eq:sect 3 3}) is
\[
\mathbf{J}_N=\left[
\begin{array}{cc}
  b^T \partial F/\partial p_l & b^T \partial F/\partial q_l \\
  b^T \partial G/\partial p_l & b^T \partial G/\partial q_l
\end{array}%
\right].
\]
Therefore, in order to preserve the underlying equilibria structure
for the numerical scheme, we need to study the matrix structure of
$\mathbf{J}_N(\tilde{p},\tilde{q})$ with $(\tilde{p},\tilde{q})\in
\mathscr{E}_H$, as what we have done before. Obviously, such arguments
apply equally to SPRK methods, and to Hamiltonian systems having
the other two kinds of equilibria structures.  In the sequel, as an
application, we give an example. Consider the following separable
1-d.o.f. Hamiltonian system
\begin{equation}\label{eq:sect 3 4}
\begin{cases}
& \dot{p}=\alpha p(1-p),\\
& \dot{q}=\alpha (2 p-1)q,
\end{cases}
\end{equation}
where $\alpha>0$. This problem has been studied in \cite{HonLiuSun},
and the solutions $p(t)\uparrow 1 (t\rightarrow +\infty)$ if
$0<p(t_0)<1$, and $p(t)\downarrow 1 (t\rightarrow +\infty)$ if
$p(t_0)>1$; while for any $q(t_0)>0$, $q(t)\uparrow +\infty
(t\rightarrow +\infty)$ if $p(t_0)>1/2$, and if $0<p_0<1/2$, $q(t)$
is first monotonically decreasing for $t<t_0=(\ln(1-p(t_0))-\ln
p(t_0)/\ln \alpha)$ and then $q(t)\uparrow +\infty$ for
$t(>t_0)\rightarrow +\infty$. It can be seen that $(0,0)$ and
$(1,0)$ are two hyperbolic equilibrium points for (\ref{eq:sect 3
4}). The stability matrices here are
\begin{equation}\label{eq:matrices}
\mathbf{J}_S(0,0)=\left[%
\begin{array}{cc}
  \alpha & 0 \\
  0 & -\alpha \\
\end{array}%
\right],\quad \mbox{and}\quad \mathbf{J}_S(1,0)=\left[%
\begin{array}{cc}
  -\alpha & 0 \\
  0 & \alpha \\
\end{array}%
\right].
\end{equation}
The symplectic Euler method for (\ref{eq:sect 3 4}) is read as
\[
\begin{cases}
& p_{l+1}=(1+z) p_l- z p_l^2,\\
& q_{l+1}=q_l+z (2 p_l-1) q_{l+1}, \quad z=\alpha h.
\end{cases}
\]
or, equivalently,
\begin{equation}\label{eq:sect 3 5}
\begin{cases}
& \delta_t^+ p_l=\alpha p_l (1-p_l),\\
& \delta_t^+ q_l=\alpha \displaystyle{\frac{(2 p_l-1)q_l}{1+z-2 z
p_l}}.
\end{cases}
\end{equation}

We have the following theorem for scheme~(\ref{eq:sect 3 5}) to give
a true simulation (see \cite{HonLiuSun})
\begin{Theorem}
The sequence $\{p_l\}\rightarrow 1$ and $\{q_l\}\rightarrow +\infty$
as $n\rightarrow +\infty$ if{f}
\begin{equation}\label{eq:condition}
0<z<1 \quad \mbox{and}\quad p(t_0)<\frac{1+z}{2z}.
\end{equation}
\end{Theorem}
If condition~(\ref{eq:condition}) is destroyed, then it is shown in
\cite{HonLiuSun} that spurious solutions or even periodic solutions
will be encountered. Through straightforward calculations, the
stability matrices of the discrete dynamic~(\ref{eq:sect 3 5}) at
its equilibrium points are given by
\[
\mathbf{J}_N(0,0)=\left[%
\begin{array}{cc}
  \alpha & 0 \\
  0 & -{\alpha}/(1+z) \\
\end{array}%
\right]\quad \mbox{and}\quad \mathbf{J}_N(1,0)=\left[%
\begin{array}{cc}
  -\alpha & 0 \\
  0 & \alpha/(1-z) \\
\end{array}%
\right].
\]
In comparison with the stability matrices in (\ref{eq:matrices}), we
see that for scheme~(\ref{eq:sect 3 5}) to preserve the hyperbolic
equilibria structure, we need require $1-z>0$, i.e., $0<z<1$.
That is, the preservation of equilibria structures only endows the
numerical integrator a prerequisite for successful simulations.
%Based on such observations, we say that the study of equilibria
%structures provides a new notion of stability for symplectic
%methods.

For Hamiltonian systems of higher dimensions, the situation will
become much more complicated. However, we can make use of the usual
linearization techniques for our investigations, which are interesting topics
for our future study.

\section*{Acknowledgements}

The work of Xiaohua Ding is supported by the NSF of China (No.
10271036) and the NSF of HIT (Weihai) (No.200518). The work of Hongyu Liu
is partly supported by NSF grant, FRG  DMS 0554571. The work of
Zaijiu Shang is supported by the NSF of China (No. 1057 1173). The
work of Geng Sun is supported by the NSF of China (No. 10471145).

\end{document}